\pdfoutput=1
\documentclass[12pt]{amsart} 

\def\htv{\widehat{\vv}}

\def\oll{\overline{\l}}							
\def\ub{\underbrace}	
\def\disp{\displaystyle}

\def\Nalr[#1,#2]{\Na_{\oll_{\rm r}(#1,#2)}(z)}
\def\Nals[#1,#2]{\Na_{\oll_\rs(#1,#2)}(z)}
\def\Nalr[#1]{\Na_{\oll_{\rm r}(#1)}(z)}
\def\Nals[#1]{\Na_{\oll_\rs(#1)}(z)}

\def\dm{\stackrel {\De} {\leftrightarrow}}
\def\dq {\stackrel {\De} {\sim}} 
\def\sdm{\stackrel {{\rm s}\De} {\leftrightarrow}}
\def\sdq {\stackrel {{\rm s}\De} {\sim}} 
\def\cyq {\stackrel{cyc}{\sim}}

\def\bgnEq{\begin{equation}}	\def\endEq{\end{equation}}
\usepackage{graphics,longtable,caption} 
\usepackage{graphicx, color, mystyle}
\usepackage{bm}

\begin{document}



\title[The Conway polynomials and Self Delta-equivalence of pretzel links]
	{The Conway polynomials and Self Delta-equivalence of pretzel links}

\author{Yasutaka NAKANISHI}

\author{Tetsuo SHIBUYA}

\author{Tatsuya Tsukamoto}

\maketitle
\begin{abstract}
	In this paper, we study the self $\De$-equivalence of pretzel links.
	If the number of	 components is $2$, then we know the complete invariants 
	in terms of the Conway polynomial for classification. We calculate the values.
	For pretzel links with $\m$ $(\geq 3)$ components, we give a necessary and 
	sufficient condition to be self $\De$-equivalent.
\end{abstract}
%
%
%
%
%
%
%
%
%
%
%
%
%
%
%
%
%
%
%
%
\section{Introduction}\label{sec1} 

	Throughout this paper, a link means a set of tame, ordered, oriented, and simple closed curves 
	in the oriented $3$-space $\bR^3$. If the number of components of a link $\k$ is $\m$, we say 
	that $\k$ is a $\m$-{\it component link}. A $1$-component link is sometimes called a knot, simply. 
	For two links $\k$ and $\l$, if there is an ambient isotopy which maps $\k$ to $\l$ preserving the 
	order and orientation, we say that $\k$ is equivalent to $\l$, denoted by $\k \approx \l$. 
 \pva
 	For two link diagrams $K$ and $L$ which differ only in one place as illustrated in Figure \ref{fig:delta}, 
	a local move between $K$ and $L$ is called a {\it $\De$-move\/}, denoted by $K \dm L$. 
	Then, for two links $\k$ and $\l$ with diagrams $K$ and $L$, $\k$ and $\l$ are said to be 
	{\it transformed into each other by a $\De$-move}. 
	Furthermore, for two links $\k$ and $\l$, if $\k$ can be transformed into $\l$ by a finite sequence of 
	$\De$-moves, $\k$ and $\l$ are said to be {\it $\De$-equivalent}, denoted by $\k \dq \l$. 
	In the case that all three arcs illustrated in Figure \ref{fig:delta} are contained in the same component, 
	the above move is called a {\it self $\De$-move\/}, denoted by $K \sdm L$. 
	Then, for two links $\k$ and $\l$ with diagrams $K$ and $L$, $\k$ and $\l$ are said to be 
	{\it transformed into each other by a self $\De$-move} \cite{S}. 
	Furthermore, for two links $\k$ and $\l$, if $\k$ can be transformed into $\l$ by a finite sequence 
	of self $\De$-moves, $\k$ and $\l$ are said to be {\it self $\De$-equivalent} (or $\De$ {\it link homotopic} 
	\cite{N, NO}), denoted by $\k \sdq \l$.
%
\bgnF	
		\iclg{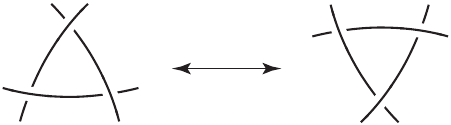}		\caption{$\De$-move}		\label{fig:delta}
\endF
%
\pva
	It is known the following result by Matveev \cite{M} and by Murakami and the first author \cite{MN}.
 
 \bgnP	
 		\label{prp:delta}
		A pair of knots $($or links$)$ are $\De$-equivalent if and only if they have the same number of 
		components and the same linking numbers between the corresponding components.
\endP
%
\bgnM
\label	{clm:sde}
		Let two links $\k$ and $\l$ be  self $\De$-equivalent. Then we have{\rm :}

 {\rm (1)} 	Each corresponding pair of sublinks of $\k$ and $\l$ are self $\De$-equivalent.
 
 {\rm (2)} 	Each corresponding pair of $2$-component sublinks  of $\k$ and $\l$
 		have the same linking numbers.
 \endM
%

	The first author and Y.Ohyama studied the self $\De$-equivalence of links in terms of 
	the coefficients of the Conway polynomial (\cite{N}, \cite{NO}). The Conway polynomial 
	$\Na_\k(z) = a_0(\k) + a_1(\k) z +$ $a_2(\k) z^2 +\cdots + a_m(\k) z^m$ is an invariant of 
	a link $\k$ which is a polynomial in a variable $z$ with integral coefficients $a_0(\k), a_1(\k), 
	a_2(\k), \cdots, a_m(\k)$. It is known that if $\k$ is a $\m$-component link $(\m \geq 2)$ 
	then $a_0(\k)  = \cdots = a_{\m-2}(\k)  = 0$, and if $i$ and $\m$ have the same parity then 
	$a_i(\k)  = 0$ (cf. \cite{K}). 
	 
\bgnP
		{\rm(\cite[Theorem~2]{N})}  \label{prp:Nak}
		If two $\m$-component links $\k = \k_1 \cup \cdots \cup \k_\m$
 		and $\l = \l_1 \cup \cdots \cup \l_\m$ are self $\De$-equivalent, then we have
 
{\rm (1)} 	$a_{\m-1}(\k) = a_{\m-1}(\l)$ {\ } and

{\rm (2)} 	$a_{\m+1}(\k) - a_{\m-1}(\k)\times \lpa \ovs{\m}{\uds{i=1}{\sum}} a_2 (\k_i) \rpa =
 		  a_{\m+1}(\l)  - a_{\m-1}(\l) \times \lpa \ovs{\m}{\uds{i=1}{\sum}} a_2 (\l_i)  \rpa $.
\endP
  
	The above is a necessary condition, but not a sufficient condition in general. For instance, 
	consider the following two $3$-component links; the $3$-component trivial link, and the split 
	sum of a $2$-component Hopf link and the unknot. They have the same Conway polynomial 
	$\Na(z) =0$. However they are not self $\De$-equivalent, since they have different linking 
	numbers of components and a (self) $\De$-move does not change the linking numbers of 
	components. But for $2$-component links, the above is also a sufficient condition as follows.
 
\bgnP
		{\rm(\cite[Theorem~3]{NO})} \label{prp:NO}
		Two $2$-component links $\k = \k_1 \cup \k_2$ and $\l = \l_1 \cup  \l_2$ are self 
		$\De$-equivalent if and only if
 
{\rm (1)} 	$a_{1}(\k) = a_{1}(\l)$  {\ } and
 
{\rm (2)} 	$a_{3}(\k) - a_{1}(\k)\times \lpa a_2 (\k_1) + a_2(\k_2) \rpa =
 		  a_{3}(\l)  - a_{1}(\l) \times \lpa a_2 (\l_1)   + a_2(\l_2) \rpa .$
\endP
 
 	In this paper, we study a necessary and sufficient condition for pretzel links to be self $\De$-equivalent. 
	For a sequence $\vv=(k_1, k_2, \ldots,  k_u)$ of non-zero integers, the unoriented link diagram 
	shown in Figure \ref{pretzel} is denoted by $P(k_1, k_2, \ldots, k_u)$ or $P_{\vv}$, which is obtained 
	as follows: First take $u$ pairs of parallel strands, introducing $k_i$ half-twists on the $i$-th pair, 
	with the convention $k_i > 0$ for right-handed half-twists and $k_i < 0$ for left-handed half-twists. 
	Then, connect each pair of parallel strands with $k_i$ and with $k_{i+1}$ by a pair of a {\it top bridge} 
	$t_i$ and a {\it bottom bridge} $b_i$, where $i=1$, $\dots$, $u$ and $k_{u+1}=k_1$. 
	We call $\vv$ a {\it pretzel sequence}.	
%
\bgnF
	\iclg{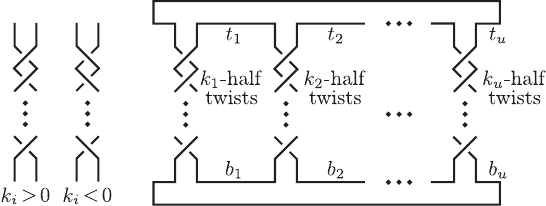}	\caption{a pretzel link diagram $P(k_1, k_2, \ldots, k_u)$}	\label{pretzel}
\endF
%
\bgnM
		\label{clm:comp}
		For a pretzel link $\k$ with a diagram $P(k_1, k_2, \ldots, k_u)$, the following can be seen.
\par
{\rm (1)} 	$\k$ is a knot if and only if all the parameters $k_i$ are odd and $u$ is odd, or exactly one 
		parameter is even. \par
{\rm (2)} 	$\k$ is a $2$-component link if and only if all the parameters $k_i$ are odd and $u$ is even, 
		or exactly two parameters are even.\par
{\rm (3)} 	$\k$ is a $\m$-component link with $\m\geq 3$ if and only if exactly $\m$ parameters are  even. 
\endM
	
	If $P(k_1, k_2, \ldots,  k_u)$ is oriented, then there are four types of orientations on each pair of 
	parallel strands as illustrated in the left four of Figure \ref{sandr}. An {\it enhanced pretzel sequence} 
	$\htv=(k_1\ep_1,$ $k_2\ep_2,$ $\ldots, k_u\ep_u)$ is a pretzel sequence with an information such 
	that $\ep_i$ is $\rs$ (resp. $\rr$) if the $i$-th pair of parallel strands with $k_i$ has anti-parallel 
	(resp. parallel) orientation ($i=1$, $2$, $\dots$, $u$). We also define $0\rs$, $0\rr$, $\infty\rs$, and 
	$\infty\rr$ as illustrated in the right eight of Figure \ref{sandr}, which we use in the following sections.
	Then our diagram is denoted by $P(k_1\ep_1, k_2\ep_2, \ldots, k_u\ep_u)$ or by $P_{\htv}$.  
	Note that as we can see in Figure \ref{sandr}, there are two possibilities to decide the orientation 
	of parallel strands with $k_i\ep_i$. Also note that if we decide an orientation of a certain top bridge 
	$t_i$, say $t_1$, then the orientations of the diagram is decided by $\htv$. Here we define the 
	orientation of $t_i$ and $b_i$ as {\it positive} (resp. {\it negative}) if the orientation is from the 
	strands with $k_i$ (resp. $k_{i+1}$) to the strands with $k_{i+1}$ (resp. $k_i$) $(i=1$, $\dots$, $u)$. 
	Then we define an oriented pretzel link diagram $P(k_1\ep_1, k_2\ep_2, \ldots, k_u\ep_u)$, 
	if one exists, as the oriented diagram such that the orientation of $t_1$ is positive.

\bgnP
	\label{prp:ori}
	There exists an oriented pretzel link diagram $P_{\htv}=P(k_1\ep_1, k_2\ep_2,$ $\ldots,$ $k_u\ep_u)$ 
	if and only if the following conditions hold:
\par
{\rm (1)} 	$\ep_i$'s are all $\rs$ 			if all the parameters $k_i$ are odd and $u$ is odd.\par
{\rm (2)} 	$\ep_i$'s are all $\rs$ or all $\rr$ 	if all the parameters $k_i$ are odd and $u$ is even.\par
{\rm (3)} 	the number of $\rr$ in $\ep_1$, $\ep_2$, $\dots$, $\ep_u$ is even, and $\ep_i$ is $\rr$ 
		when $k_i$ is odd if there is an even integer in $k_1, \dots, k_u$.
\endP
	Note that by Proposition \ref{prp:ori}, if there exists an oriented pretzel link diagram $P(k_1\ep_1, 
	k_2\ep_2,$ $\ldots,$ $k_u\ep_u)$, then the number of $\rr$ in $\ep_1$, $\ep_2$, $\dots$, $\ep_u$ is even.	

%
\bgnF
	\iclg{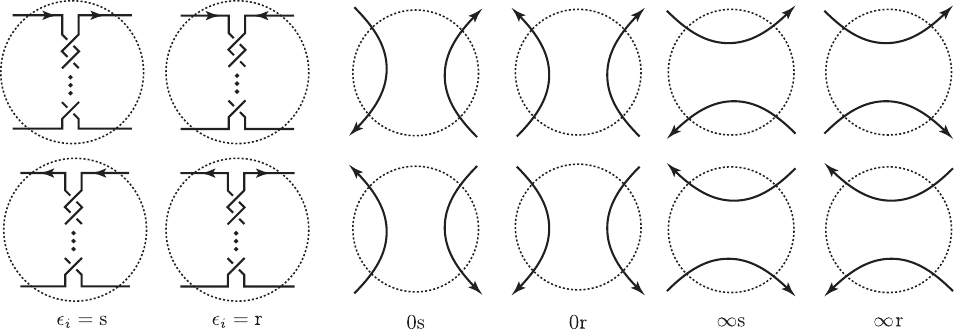}	\caption{anti-parallel and parallel strands}		\label{sandr}
\endF
	
	Moreover, define an equivalence relation $\cyq$ on the set of enhanced pretzel 
	sequences as $(k_1\ep_1, \dots, k_u\ep_u)$ $\cyq$ $(\ell_1\ep'_1, \dots, \ell_u\ep'_u)$ 
	if and only if there is a number $t$ $(1\leq t\leq u)$ such that $\ell_i\e'_i=k_{i+t}\ep_{i+t}$ 
	$(k_{u+1}\ep_{u+1}=k_1\ep_1)$ or $\ell_i\ep'_i=k_{u-i+t}\ep_{u-i+t}$ 	for any $i$ 
	$(i= 1, 2, \dots, u)$. Then, we have the following.

\bgnL
		\label{lem:EPS} 
		Let $\htv$ and $\htv'$ be two enhanced pretzel sequences such that $\htv$ $\cyq$ $\htv'$. 
		Then we have that $P_{\htv}$ $\approx$ $P_{\htv'}$.
\endLR
		The statement holds, since we have the following by considering the diagrams:
\pva
		$P(k_1\ep_1, k_2\ep_2, \dots, k_u\ep_u)$ $\approx$ $P(k_2\ep_2, \dots, k_u\ep_u, k_1\ep_1)$ and \pva
		$P(k_1\ep_1, k_2\ep_2, \dots, k_u\ep_u)$ $\approx$ $P(k_u\ep_u, \dots, k_2\ep_2, k_1\ep_1)$.
\endR

	We also define $\r(\htv)$ as $\uds{k_i: {\rm odd}}{\sum} k_i - \sharp \{ i \vert k_i = -2\}$.
	For an enhanced pretzel sequence $\htv=(k_1\ep_1, k_2\ep_2, \dots, k_u\ep_u)$ with even 
	parameters, let $\htv_e=(k_{i_1}\ep_{i_1},$ $k_{i_2}\ep_{i_2}$, $\dots$, $k_{i_\m}\ep_{i_\m})$ be the 
	subsequence of $\htv$ consisting of all the elements with even parameters  ($i_1<i_2<\cdots<i_\m$). 
	Note that the number of components of $P_{\htv}$ is $\m$ by Claim \ref{clm:comp} (3). 
	Now the following is our main theorem.

\bgnT
\label	{thm:EPS} 
		Let $\htv$ and $\htv'$ be two enhanced pretzel sequences.
		If $\k$ and  $\l$ are $\m$-component pretzel links $(\m\geq 3)$ with diagrams $P_{\htv}$ 
		and $P_{\htv'}$, respectively, then, $\k$ and $\l$ are self $\De$-equivalent if and only if
\pvan
{\rm (1)} 	$\htv_e$ $\cyq$ $\htv'_e$ by permitting $2\rr = -2\rs$ and $2\rs = -2\rr$, and\pvan
{\rm (2)} 	$\r(\htv)=\r(\htv')$.
 \endT	
	
	For example, we know that $P(4\rs, 5\rr, 6\rr, -2\rr, -3\rr)\sdq P(6\rr, 2\rs, 7\rr, 4\rs, -5\rr, -1\rr)$
	by Theorem \ref{thm:EPS}.
\pvc	
	The concordance of links is defined in \cite{M}. For a links $\k$ , if $\k$ is concordant to a trivial link, 
	the link is called a {\it slice link\/}. Then we have the following.

\bgnT
\label	{thm:Slice} 
		A slice pretzel link is self $\De$-equivalent to a trivial link.
\endT

	The converse of Theorem \ref{thm:Slice} is not valid. A $2$-component pretzel link $\k$ with a 
	diagram $P(2k\ep, -2k\ep, 1\rr, -3\rr, -5\rr, 7\rr)$	is self $\De$-equivalent to the $2$-component 
	trivial link by Theorem \ref{thm:EPS}.	 However, $\k$ is not slice, since $\k$ consists of the unknot 
	and a non-slice knot.
\pvc
	This paper is organized as follows. 
	In Section \ref{secdef}, we review definitions and preliminaries to show Theorem~\ref{thm:EPS}.
	In Section \ref{sec:main}, we consider self $\De$-equivalence for $\m$-component pretzel links	$(\m \geq 3)$.
	In Section \ref{sec:2comp}, we consider self $\De$-equivalence for $2$-component pretzel links.
	In Section \ref{sec:eras}, we  give a proof of Theorem \ref{thm:Slice}.
	
%
%
%
%
%
%
%
%
%
%
%
%
%
%
%
%
%
%
%
%

\section{Definitions and Preliminaries}\label{secdef} 

	
	First we prove Proposition \ref{prp:ori} and in the rest of this section, we calculate the Conway
	polynomial of links.
\bgnR
	[Proof of Proposition \ref{prp:ori}.]  
	Note that the unoriented pretzel link diagram $P_{\vv}$ admits $2^\m$ types of orientations, where
	$\m$ is the number of components of $P_{\vv}$, which is equal to the number of even parameters
	of $\vv$ if exists and $1$ if all the parameters $k_i$ are odd by Claim \ref{clm:comp}. 
\pvan
	(1) 	
	Since all the parameters $k_i$ are odd and $u$ is odd, we have that $\m=1$, and if we walk along 
	$P_{\vv}$ from $t_1$ with a positive (resp. negative) orientation, then we meet $b_2$, $t_3$, $\dots$, 
	$b_{u-1}$, $t_u$, $b_1$, $t_2$, $\dots$, $t_{u-1}$, $b_u$ in this order (resp. in the opposite order of
	this) and come back to $t_1$. 
	$(\Rightarrow)$
	Since $\m=1$, $P_{\vv}$ admits $2$ types of orientations, that is realized by assigning $t_1$ a positive
	or negative orientation. In the former (resp. latter) case, we know that the other top and bottom bridges 
	are assigned positive (resp. negative) orientations by the above observation. Thus, the strand of any 
	$k_i$ has anti-parallel orientation, i.e., each $\ep_i$ is $\rs$, since $t_{i-1}$ and $b_{i-1}$ have the 
	same orientations $(i=1,\dots, u)$.
	$(\Leftarrow)$ 
	By assigning an orientation to $P_{\vv}$ so that $t_1$ has a positive orientation, we know that the oriented
	diagram is $P_{\htv}$.
\pvan
(2) 	Since all the parameters $k_i$ are odd and $u$ is even, we have that $\m=2$, and we know the following.
	If we walk along $P_{\vv}$ from $t_1$ with a positive (resp. negative) orientation, then we meet $b_2$, 
	$t_3$, $\dots$, $t_{u-1}$, $b_u$ in this order (resp. in the opposite order of this) and come back to $t_1$. 
	If we walk along $P_{\vv}$ from $b_1$ with a positive (resp. negative) orientation, then we meet $t_2$, 
	$b_3$, $\dots$, $b_{u-1}$, $t_u$ in this order (resp. in the opposite order of this) and come back to $b_1$. 
	$(\Rightarrow)$ 
	Since $\m=2$, $P_{\vv}$ admits $2^2$ types of orientations that is realized by assigning $t_1$ and $b_1$
	the same or opposite orientations. In the former (resp. latter) case, the strand of any $k_i$ has anti-parallel 
	(resp. parallel) orientation, i.e., each $\ep_i$ is $\rs$ (resp. $\rr$), since $t_{i-1}$ and $b_{i-1}$ have the 
	same (resp. opposite) orientations $(i=1,\dots, u)$. Hence the condition holds.
	$(\Leftarrow)$ 
	If $\ep_i$'s are all $\rs$ (resp. $\rr$), then assign an orientation to $P_{\vv}$ so that $t_1$ has a positive 
	orientation and $b_1$ has a positive (resp. negative) orientation. We know that the oriented diagram is $P_{\htv}$.
\pvan
(3)
	Let $k_{i_1}$, $k_{i_2}$, $\dots$, $k_{i_\m}$ be the even parameters of $\vv$ and $K_1$, $K_2$, $\dots$,
	$K_\m$ the knot components of $P_{\vv}$ such that the crossings of $K_j$ are of $k_{i_j}$, $\dots$, $k_{i_{j+1}-1}$	
	$(j=1,\dots, \m)$. If $i_{j+1}-i_j$ is odd and we walk along $K_j$ from $t_{i_j}$ with a positive (resp. negative) 
	orientation, then we meet $b_{i_j+1}$, $t_{i_j+2}$, $\dots$, $b_{i_{j+1}-2}$, $t_{i_{j+1}-1}$, $b_{i_{j+1}-1}$, 
	$t_{i_{j+1}-2}$, $\dots$, $t_{i_j+1}$, $b_{i_j}$ in this order 	(resp. in the opposite order of this) and come back to 
	$t_{i_j}$. If $i_{j+1}-i_j$ is even and we walk along $K_j$ from $t_{i_j}$ with a positive (resp. negative) orientation, 
	then we meet $b_{i_j+1}$, $t_{i_j+2}$, $\dots$, 	$t_{i_{j+1}-2}$, $b_{i_{j+1}-1}$, $t_{i_{j+1}-1}$, $b_{i_{j+1}-2}$, 
	$\dots$, $t_{i_j+1}$, $b_{i_j}$ in this order (resp. in the opposite order of this) and come back to $t_{i_j}$. 
	Note that in any of the above four cases, each $\ep_h$ is $\rr$ $(i_j+1\leq h\leq i_{j+1}-1)$.
\pvan
	$(\Rightarrow)$	
	Assume that $P_{\htv}$ is realized and take any odd parameter $k_a$. Let $K_\a$ be the knot component of 
	$P_{\htv}$ which contains the strands of $k_a$. Then we know that $\ep_h$ is $\rr$ $(i_\a+1\leq h\leq i_{\a+1}-1)$ 
	by the above observation, and thus that $\ep_a$ is $\rr$. Now note that the orientations of $t_{i-1}$ and $t_i$ 
	are the same and the opposite if $\ep_i$ is $\rs$ and $\rr$, respectively. Thus the orientation of $t_{\m+1}$ is 
	the same (resp. opposite) as that of $t_1$ if the number of $\rr$ in $\ep_1$, $\ep_2$, $\dots$, $\ep_u$ is even 
	(resp. odd). Therefore, the number of $\rr$ in $\ep_1$, $\ep_2$, $\dots$, $\ep_u$ is even, since $t_{\m+1}=t_1$.
\pvan
	$(\Leftarrow)$ 
	We may assume that $k_{i_1}=k_1$. We assign $P_{\vv}$ an orientation to realize $\htv$. First assign $K_1$ 
	an orientation by assigning $t_1$ a positive orientation. Then we have that $t_{i_2-1}$ is assigned a positive 
	(resp. negative) orientation if $i_2-i_1=i_2-1$ is odd (resp. even). Assume that $K_1$, $\dots$, $K_{\a-1}$ are 
	oriented $(2\leq \a\leq u)$. Then assign $K_\a$ an orientation by assigning $t_{i_\a}$ the same (resp. opposite) 
	orientation as that of $t_{i_\a-1}$ if $\ep_{i_\a}$ is $\rs$ (resp. $\rr$). Finally assign $K_{\m}$ an orientation by 
	assigning $t_{i_\m}$ the same (resp. opposite) orientation as that of $t_{i_\m-1}$ if $\ep_{i_\m}$ is $\rs$ (resp. 
	$\rr$). Note that the orientations of $t_{i-1}$ and $t_i$ are the same and the opposite if $\ep_i$ is $\rs$ and $\rr$, 
	respectively. Since the number of $\rr$ in $\ep_1$, $\ep_2$, $\dots$, $\ep_u$ is even, the orientation of 
	$t_{\m+1}(=t_1)$ is positive, and thus the orientation of $P_{\htv}$ is well defined. 
\endR
%
%
%
%
%
%
%
%
%
%
%
%
%
%
%
%
%
%
%
%
	\subsection{the Conway polynomial of links}\label{sbscwy} 
	
	Let $\l(k\ep)$ be a link which has an oriented diagram $L(k\ep)$ containing a strand $k\ep$ of Figure
	\ref{sandr} ($\ep=\rs$ or $\rr$). Moreover let $\l(k_1\ep_1, k_2\ep_2, \dots, k_n\ep_n)$ $(n\geq 2)$ be 
	a link which has an oriented diagram $L(k_1\ep_1, k_2\ep_2, \dots, k_n\ep_n)$ obtained by connecting 
	the right top (resp. bottom) of the strand for $k_i\ep_i$ and left top (resp. bottom) of the strand for 
	$k_{i+1}\ep_{i+1}$ $(i=1, \dots, n-1)$, where each $k_i\ep_i$ is of Figure \ref{sandr}. Here we retake the 
	broken circle to surround $L(k_1\ep_1, k_2\ep_2, \dots, k_n\ep_n)$ (see Figure \ref{sum} for an example).
	
\bgnF
	\iclg{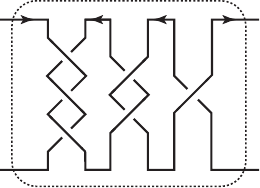}		\caption{$L(-3\rr, 2\rs, 1\rr)$}	\label{sum}
\endF

	If these links are in a common equation, then we assume that these links have corresponding oriented 
	diagrams which are identical outside the broken circles.
\pvc	
	The Conway polynomial $\Na_\k(z)$ of a link $\k$ can be determined axiomatically by the following two 
	equations (cf. \cite{K}):\pvan
	(1) 	For the unknot $O$, $\Na_O(z) = 1$ and \pvan	
	(2) 	$\Na_{\l(1\rr)}(z)-\Na_{\l(-1\rr)}(z)=z\Na_{\l(0\rr)}(z)$.
\pvc
	In the following, we may denote the Conway polynomial of $\k$ simply by $\Na_\k$.
	Let $\phi_i(t)=\pY t+i\\ 2i+1 \epy$ and $\psi_i(t)=\pY t+i\\ 2i \epy$, where\pva
\bgnC	
	$\pY \a\\ ~n\epy=\Frc{\a(\a-1)\cdots(\a-n+1)}{n!}$ ($\a\in\bR$, $n\in\bN$) and $\pY \a\\ ~0\epy=1$.
\endC	

	Note that the following equation holds.
\bgnEq
	\label{eqPascal} \pY \a\\ ~n+1\epy	+\pY \a\\ ~n\epy =\pY \a+1\\ ~n+1\epy.	\endEq	
	
	Then, define polynomials $\Phi(t,z)$ and $\Psi(t,z)$ in variable $z$ for an integer $t$ as
\bgnC	
	$\Phi(t,z)=\Sum_{i=0}^{\infty} \phi_i(t) z^{2i+1}$ and 
	$\Psi(t,z)=\Sum_{i=0}^{\infty} \psi_i(t) z^{2i}$.
\endC
	Here note that $\phi_i(t)=0$ if and only if $i\geq |t|$ and $\psi_i(t)=0$ if and only if
	$i\geq t+1$ if $t\geq 0$ and $i\geq -t$ if $t<0$.
	
\bgnL
\label	{lem:Pascal}		The following holds.\pvan
$(1)$ 	$\Phi(-p,z)=-\Phi(p,z)$, \pvan
$(2)$ 	$\Psi(-p-1,z)=\Psi(p,z)$, \pvan
$(3)$		$\Phi(p,z)+z\Psi(p,z)=\Phi(p+1,z)$, and		\pvan
$(4)$		$\Psi(p-1,z)+z\Phi(p,z)=\Psi(p,z)$
\endLR
	(1) (2) Both equations hold, since we have that $\phi_i(-p)=-\phi_i(p)$ and
		$\psi_i(-p-1,z)=\psi_i(p,z)$ by definition. 
	(3) (4) Both equations are obtained by Equation (\ref{eqPascal}) as follows.\pvan
	(3) 	$\ovs{\infty}{\uds{i=0}{\sum}} \py p+i\\ 2i+1 \epy\! z^{2i+1}+
		z\ovs{\infty}{\uds{i=0}{\sum}} \py p+i\\ 2i	 \epy\! z^{2i}=
		\ovs{\infty}{\uds{i=0}{\sum}} \lbc \py p+i\\ 2i+1 \epy\!+\!\py p+i\\ 2i \epy\rbc\! z^{2i+1}
		=\ovs{\infty}{\uds{i=0}{\sum}}  \py p+1+i\\ 2i+1  \epy\! z^{2i+1}$.\pvcn
	(4)	$\ovs{\infty}{\uds{i=0}{\sum}} \py p-1+i\\ 2i \epy\! z^{2i}+
		z\ovs{\infty}{\uds{i=0}{\sum}} \py p+i\\ 2i+1 \epy\! z^{2i+1}=
		\py p-1\\~0 \epy z^0+\ovs{\infty}{\uds{i=1}{\sum}} \py p-1+i\\ 2i \epy\! z^{2i}+
		\ovs{\infty}{\uds{i=1}{\sum}} \py p+i-1\\ 2i-1 \epy\! z^{2i}$\pvc
		$=\py p\\~0 \epy z^0+\ovs{\infty}{\uds{i=1}{\sum}} \py p+i\\ 2i \epy\! z^{2i}
		=\ovs{\infty}{\uds{i=0}{\sum}} \py p+i\\ 2i \epy\! z^{2i}$.
\endR
\bgnL
\label	{lem:TwRed}	
		We have the following for an integer $p$. \pvan
$(1)$		$\Na_{\l(2p\rs)}	  = 	\Na_{\l(0\rs)}-pz\Na_{\l(\infty\rs)}$, \pvan
$(2)$		$\Na_{\l(2p+1\rs)}= 	\Na_{\l(1\rs)}-pz\Na_{\l(\infty\rr)}$, \pvan
$(3)$		$\Na_{\l(2p\rr)}	= 	\Phi(p,z)\Na_{\l(1\rr)}	+	\Psi(p-1,z)\Na_{\l(0\rr)}$, and \pvan
$(4)$		$\Na_{\l(2p+1\rr)}= 	\Psi(p,z)\Na_{\l(1\rr)} +	\Phi(p,z)	\Na_{\l(0\rr)}$.
\endLR
		Each equation is obtained by induction on $p$ applying the recursive formula of the
		Conway polynomial, i.e., $\Na_{\l(1\rr)}-\Na_{\l(-1\rr)}=z\Na_{\l(0\rr)}$. In fact, we use
		the following equation obtained by the formula.
\bgnEq
		\label{eqCRFs}\Na_{\l(t\pm 2\rs)}=
		\bgnseq
		\Na_{\l(t\rs)}\mp z\Na_{\l(\infty\rs)}~~\mathrm{if}~~ t~ \mathrm{: even}\\[1ex]
		\Na_{\l(t\rs)}\mp z\Na_{\l(\infty\rr)}~~\mathrm{if}~~ t ~\mathrm{: odd}  \endseq\endEq	
\bgnEq
		\label{eqCRFr}\Na_{\l(t\pm 2\rr)}=\Na_{\l(t\rr)}\pm z\Na_{\l(t\pm1\rr)}		\endEq			

	(1)(2) The equations for $p=0$ clearly hold. If $p>0$ (resp. $p<0$) and the equations for $p$ hold, 
		then we know that the equations for $p+1$ (resp. $p-1$) hold by Equation (\ref{eqCRFs}). \pva
	(3)(4) The equations for $p=0$ hold, since $\Phi(0,z)=0$, $\Psi(-1,z)=1$, and $\Psi(0,z)=1$. 
		
		Then we have the following by Equation (\ref{eqCRFr}) and Lemma \ref{lem:Pascal}.\pvan
		$p>0:$
		$\Na_{\l(2(p+1)\rr)}=\Na_{\l(2p\rr)}+ z\Na_{\l(2p+1\rr)}$	\pvb\hsC
		$=\lpa\Phi(p,z)+z\Psi(p,z)\rpa\Na_{\l(1\rr)}+\lpa\Psi(p-1,z)+z\Phi(p,z)\rpa\Na_{\l(0\rr)}$\pvb\hsC
		$=\Phi(p+1,z)\Na_{\l(1\rr)}+\Psi(p,z)\Na_{\l(0\rr)}$.
\pbg
		$\Na_{\l(2(p+1)+1\rr)}=\Na_{\l(2p+1\rr)}+ z\Na_{\l(2(p+1)\rr)}$	\pvb\hsC
		$=\lpa\Psi(p,z)+z\Phi(p+1,z)\rpa\Na_{\l(1\rr)}+\lpa\Phi(p,z)+z\Psi(p,z)\rpa\Na_{\l(0\rr)}$\pvb\hsC
		$=\Psi(p+1,z)\Na_{\l(1\rr)}+\Phi(p+1,z)\Na_{\l(0\rr)}$.
\pvbn
		$p<0:$
		$\Na_{\l(2(p-1)\rr)}=\Na_{\l(2p\rr)}- z\Na_{\l(2p-1\rr)}$	\pvb\hsC
		$=\lpa\Phi(p,z)-z\Psi(p-1,z)\rpa\Na_{\l(1\rr)}+\lpa\Psi(p-1,z)-z\Phi(p-1,z)\rpa\Na_{\l(0\rr)}$\pvb\hsC
		$=\Phi(p-1,z)\Na_{\l(1\rr)}+\Psi(p-2,z)\Na_{\l(0\rr)}$.
\pbg
		$\Na_{\l(2(p-1)+1\rr)}=\Na_{\l(2p+1\rr)}- z\Na_{\l(2p\rr)}$	\pvb\hsC
		$=\lpa\Psi(p,z)-z\Phi(p,z)\rpa\Na_{\l(1\rr)}+\lpa\Phi(p,z)-z\Psi(p-1,z)\rpa\Na_{\l(0\rr)}$\pvb\hsC
		$=\Psi(p-1,z)\Na_{\l(1\rr)}+\Phi(p-1,z)\Na_{\l(0\rr)}$.
\endR	


	By connecting the left (resp. right) top and bottom of a diagram of Figure \ref{sandr} by a simply 
	connected arc if the orientation is realizable, we obtain a link diagram $T(k\ep)$ of the $(2,k)$-torus 
	link $\tau(k\ep)$. Since $\tau(1\rr)$ and $\tau(\infty\rs)$ (resp. $\tau(0\rr)$ and $\tau(0\rs)$) are the trivial 
	knots (resp. the trivial $2$-component links), and hence $\Na_{\tau(1\rr)}(z)=\Na_{\tau(\infty\rs)}(z)=1$ 
	(resp. $\Na_{\tau(0\rr)}(z)=\Na_{\tau(0\rs)}(z)=0$), we have the following by Lemma \ref{lem:TwRed}.

\bgnO
\label	{cor:TorusAlex} The following holds for $\Na_{\tau(k{\ep})}$.\pva
$(1)$		$\Na_{\tau(2p\rs)} =-pz$, \par
$(2)$		$\Na_{\tau(2p\rr)} = \Phi(p,z)=pz+\Frc{1}{\,6\,}p(p^2-1)z^3+\cdots$, and  \par
$(3)$		$\Na_{\tau(2p+1\rr)} = \Psi(p,z)=1+\Frc{1}{\,2\,}p(p+1)z^2+\cdots$.
\endO
	Hence, we know that $(2,2p)$-torus links possess mutually distinct Conway polynomials
	except the cases that $\Na_{\tau(0\rs)}=0=\Na_{\tau(0\rr)}$, $\Na_{\tau(\pm 2\rs)}=\mp z=\Na_{\tau(\mp2\rr)}$.	
	In fact, we have that $\tau(0\rs) \approx \tau(0\rr)$ and that $\tau(\pm 2\rs) \approx \tau(\mp 2\rr)$.
	Moreover we have the following by Proposition \ref{prp:NO} and Corollary \ref{cor:TorusAlex}.
%
\bgnO
\label	{cor:TorusDelta}
		Links with diagrams $\tau(2p\ep)$ and $\tau(2q\ep')$, where $p\neq q$ or $\ep\neq\ep'$
		$(p,q = 0, \pm 1, \pm 2, \dots; \ep,\ep' = \rs, \rr)$ are not self $\De$-equivalent 
		except $\tau(0\rs) \approx \tau(0\rr)$, $\tau(\pm 2\rs) \approx \tau(\mp 2\rr)$.
\endO


	Then we have the following by Lemmas \ref{lem:Pascal} and \ref{lem:TwRed}. We note that equation 
	(4) has been obtained in terms of the Alexander polynomial in \cite{STUI} (Theorem 2.4).
	
\bgnP
\label	{prp:ErsCwy}
		We have the following for an integer $p$. \pvan
$(1)$		$\Na_{\l(2p\rs, -2p\rs)}=\Na_{\tau(2p\rs)}\Na_{\tau(-2p\rs)}\Na_{\l(\infty\rs)}=-p^2z^2\Na_{\l(\infty\rs)}$, \pvan
$(2)$		$\Na_{\l(2p+1\rs, -(2p+1)\rs)} = 	(1+p^2z^2)\Na_{\l(\infty\rr)}$, \pvan
$(3)$		$\Na_{\l(2p\rr, -2p\rr)}=\Na_{\tau(2p\rr)}\Na_{\tau(-2p\rr)}\Na_{\l(\infty\rs)}= -\Phi^2(p,z)\Na_{\l(\infty\rs)}$, and \pvan
$(4)$		$\Na_{\l(2p+1\rr, -(2p+1)\rr)}=\Na_{\tau(2p+1\rr)}\Na_{\tau(-(2p+1)\rr)}\Na_{\l(\infty\rs)}= \Psi^2(p,z)\Na_{\l(\infty\rs)}$.
\endPR
	(1)	Note that $L(\infty\rs,\infty\rs)=L(\infty\rs)$, that $L(0\rs, 0\rs)$ is a split diagram, and thus
		that $\Na_{\l(0\rs, 0\rs)}=0$. Moreover, note that $\l(0\rs, \infty\rs)\approx\l(\infty\rs, 0\rs)$.
		The equation is obtained by applying Lemma \ref{lem:TwRed} (1) twice.\pva
		$\Na_{\l(2p\rs, -2p\rs)}=\Na_{\l(0\rs, -2p\rs)}-pz\Na_{\l(\infty\rs, -2p\rs)}$\pvb\hsB
		$=\Na_{\l(0\rs, 0\rs)}+pz(\Na_{\l(0\rs, \infty\rs)}-\Na_{\l(\infty\rs, 0\rs)})-p^2z^2\Na_{\l(\infty\rs, \infty\rs)}$
		$=-p^2z^2\Na_{\l(\infty\rs)}$.
\pvcn
	(2)	Note that $\l(1\rs, -1\rs)\approx\l(\infty\rr)$, $\l(1\rs, \infty\rr)\approx\l(\infty\rr, 1\rs)$, and that 
		$\l(\infty\rr, \infty\rr)\approx\l(\infty\rr)$. We have the following by Lemma \ref{lem:TwRed} (2) and 
		Equation (2.2), since $-(2p+1)=2(-p-1)+1$.\pva
		$\Na_{\l(2(-p-1)+1\rs)} = \Na_{\l(1\rs)}-(-p-1)z\Na_{\l(\infty\rr)}$\pbe\hsB
		$=\Na_{\l(-1\rs)}-z\Na_{\l(\infty\rr)}+(p+1)z\Na_{\l(\infty\rr)}$
		$=\Na_{\l(-1\rs)}+pz\Na_{\l(\infty\rr)}$.\pva
		Hence we have the following.\pva	
		$\Na_{\l(2p+1\rs, -(2p+1)\rs)} = \Na_{\l(1\rs, -(2p+1)\rs)}-pz\Na_{\l(\infty\rr, -(2p+1)\rs)}$\pvb
		$=(\Na_{\l(1\rs, -1\rs)}+pz\Na_{\l(1\rs, \infty\rs)})-pz(\Na_{\l(\infty\rr, 1\rs)}+(p+1)z\Na_{\l(\infty\rr, \infty\rr)})$
		$=(1-p(p+1)z^2)\Na_{\l(\infty\rr)}$.
\pvcn
	(3)	Note that $\l(1\rr, -1\rr)\approx\l(\infty\rs)$, that $\l(0\rr, 1\rr)\approx\l(1\rr, 0\rr)\approx\l(0\rs)$,
		and that $L(0\rr, 0\rr)$ is a split diagram, and thus that $\Na_{\l(0\rr, 0\rr)}=0$. We have the 
		following by Lemma \ref{lem:TwRed} (3), Equation (2.2), and Lemma \ref{lem:Pascal} (4), since
		$-(2p)=2(-p)$. 
\npg
		$\Na_{\l(2(-p)\rr)} = \Phi(-p,z)\Na_{\l(1\rr)}+\Psi(-p-1,z)\Na_{\l(0\rr)}$\pvb\hsB
		$= \Phi(-p,z) \Na_{\l(-1\rr)}+ (z\Phi(-p,z)+\Psi(-p-1,z))\Na_{\l(0\rr)}$\pvb\hsB
		$= \Phi(-p,z) \Na_{\l(-1\rr)}+ \Psi(-p,z)\Na_{\l(0\rr)}$.\pva
		Hence we have the following by Lemma \ref{lem:Pascal} (1)(2).\pva		
		$\Na_{\l(2p\rr,-2p\rr)}	=\Phi(p,z)\Na_{\l(1\rr,-2p\rr)}+\Psi(p-1,z)\Na_{\l(0\rr,-2p\rr)}$\pvb\hsA
		$=\Phi(p,z)(\Phi(-p,z) \Na_{\l(1\rr,-1\rr)}+ \Psi(-p,z)\Na_{\l(1\rr,0\rr)})$\pvb\hsB
			$+\Psi(p-1,z)(\Phi(-p,z)\Na_{\l(0\rr,1\rr)}+\Psi(-p-1,z)\Na_{\l(0\rr,0\rr)})$\pvb\hsA
		$=\Phi(p,z)\Phi(-p,z) \Na_{\l(\infty\rs)}+ (\Phi(p,z)\Psi(-p,z)+\Psi(p-1,z)\Phi(-p,z))\Na_{\l(0\rs)}$\pvb\hsB
		$+\Psi(p-1,z)\Psi(-p-1,z)\Na_{\l(0\rr,0\rr)}$\pvb\hsA
		$=-\Phi^2(p,z)\Na_{\l(\infty\rs)}+ (\Phi(p,z)\Psi(-p,z)-\Psi(-p,z)\Phi(p,z))\Na_{\l(0\rs)}+0$
		$=-\Phi^2(p,z)\Na_{\l(\infty\rs)}.$
\pvcn
	(4)	We have the following by Lemma \ref{lem:TwRed} (4), Equation (2.2), and Lemma 
		\ref{lem:Pascal} (3), where note that $-(2p+1)=2(-p-1)+1$. \pva
		$\Na_{\l(2(-p-1)+1\rr)} = \Psi(-p-1,z)\Na_{\l(1\rr)}+\Phi(-p-1,z)\Na_{\l(0\rr)}$\pvb\hsB
		$= \Psi(-p-1,z) \Na_{\l(-1\rr)}+ (z\Psi(-p-1,z)+\Phi(-p-1,z))\Na_{\l(0\rr)}$\pvb\hsB
		$= \Psi(-p-1,z) \Na_{\l(-1\rr)}+ \Phi(-p,z)\Na_{\l(0\rr)}$.
\pvc
		Hence we have the following by Lemma \ref{lem:Pascal} (1)(2).\pva		
		$\Na_{\l(2p+1\rr,-(2p+1)\rr)}=\Psi(p,z)\Na_{\l(1\rr,-(2p+1)\rr)}+\Phi(p,z)\Na_{\l(0\rr,-(2p+1)\rr)}$\pvb\hsA
		$=\Psi(p,z)(\Psi(-p-1,z) \Na_{\l(1\rr,-1\rr)}+ \Phi(-p,z)\Na_{\l(1\rr,0\rr)})$\pvb\hsC
			$+\Phi(p,z)(\Psi(-p-1,z)\Na_{\l(0\rr,1\rr)}+\Phi(-p-1,z)\Na_{\l(0\rr,0\rr)})$\pvb\hsA
		$=\Psi(p,z)\Psi(-p-1,z) \Na_{\l(\infty\rs)}+ (\Psi(p,z)\Phi(-p,z)$\pvb\hsC
		$+\Phi(p,z)\Psi(-p-1,z))\Na_{\l(0\rs)}+\Phi(p,z)\Phi(-p-1,z)\Na_{\l(0\rr,0\rr)}$\pvb\hsA
		$=\Psi^2(p,z)\Na_{\l(\infty\rs)}+ (-\Psi(p,z)\Phi(p,z)+\Phi(p,z)\Psi(p,z))\Na_{\l(0\rs)}+0$\pvb\hsA
		$=\Psi^2(p,z)\Na_{\l(\infty\rs)}.$
		\endR

	In the following, we may use the same symbol for a link and diagrams of it unless stated.
%
%
%
%
%
%
%
%
%
%
%
%
%
%
%
%
%
%
%
%

\section{Self Delta-equivalence of pretzel links}\label{sec:main}

%

	The following lemma is a key for the proof of Theorem \ref{thm:EPS} . Note that a pair of parallel 
	strands with odd number of half twists of a pretzel link with more than $2$ components belong
	to the same component.
\pvc
	Recall that $\l(k_1\ep_1, k_2\ep_2, \dots, k_n\ep_n)$ is a link which has a diagram 
	$L(k_1\ep_1, k_2\ep_2, \dots, k_n\ep_n)$ obtained by connecting the right top (resp. bottom) of 
	$k_i\ep_i$ and left top (resp. bottom) of $k_{i+1}\ep_{i+1}$, where each $k_i$ is of Figure \ref{sandr}.
	Then we have the following.
\bgnL
\label	{lem:key}
		For an odd integer $k$ $(|k|\!\geq\! 3)$, we have that $\l(k)\!\dq\!\l(\ub{e,\dots,e}_{\vert k \vert})$, 
		where $\disp e = \!{k \over {\vert k \vert}}$.
 \endLR
		Figure \ref{fig:vtoh} illustrates that $\Frc{1}{2}(|k|-1)$ times of $\De$-moves realize that 
		$\l(k) \dq \l(e, e, k-2e)$ for a positive integer $k$. We obtain the conclusion by applying this
		process inductively. Similarly, we can obtain the conclusion for a negative odd integer $k$.
\endR
\bgnF
\iclg		{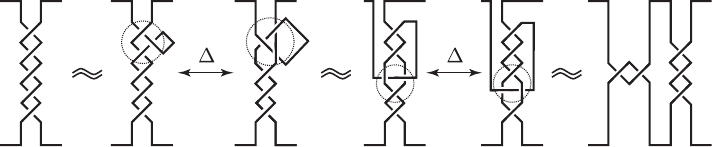}
\caption	{transformation of vertical twists into horizontal twists}
\label	{fig:vtoh}
\endF
%
\bgnM
\label	{clm:normal}
 		$\l(\cdots, -2\ep, \cdots)\, {\approx}\, \l(\cdots, 2\ep', (-1)\rr, \cdots)$
 		where $(\ep, \ep') = (\rs, \rr)$ or $(\rr, \rs)$.
 \endMR
		Figure \ref{flype1} illustrates  the required result.
\endR
\bgnF
		\iclg{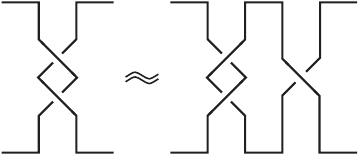}	\caption{}		\label{flype1}
\endF

 	It can be seen that an order of $k_i$ and $\pm 1$ can be switched.
 
\bgnM
\label	{clm:flype}
 		$\l(\cdots, \pm 1\rr, k\ep, \cdots) \approx
 		\l(\cdots, k\ep, \pm 1\rr, \cdots)$, where $\ep \in \{\rs, \rr\}$.
\endMR	
 		Figure \ref{flype2} illustrates  the required result.
\endR\bgnF
		\iclg{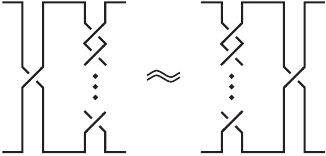}	\caption{}		\label{flype2}
\endF

	An (enhanced) pretzel sequence $(k_1\ep_1, k_2\ep_2, \dots, k_\m\ep_\m)$ is {\it standard} 
	if $k_i$ is even and is not $-2$ for any $i$ $(1\leq i\leq \m)$. For a standard enhanced pretzel 
	sequence	$\htv=(k_1\ep_1, k_2\ep_2, \dots, k_\m\ep_\m)$, let $P_{\htv}^m=P(k_1\ep_1, 
	k_2\ep_2, \dots, k_\m\ep_\m,$ $e_1\rr$, $\dots$, $e_{|m|}\rr)$, where $m$ is a non-zero 
	integer and $\disp e_i= {m \over {\vert m \vert}}$ $(1\leq i\leq \vert m \vert)$.

\bgnP\label{prp:EPS}
		Assume that $\m\geq 3$ and that enhanced pretzel sequences $\htv=(k_1\ep_1$, 
		$k_2\ep_2$, $\dots$, $k_\m\ep_\m)$ and $\htv'=( l_1\d_1$, $ l_2\d_2$, $\dots$, $ l_\m\d_\m)$ 
		are standard. Then, the following three conditions are equivalent.		\pva
{\rm (i)} 	$P_{\htv}^m\,\approx\,P_{\htv'}^{m'}$. \par\vspace{1mm}
{\rm (ii)} 	$P_{\htv}^m\,\sdq\,P_{\htv'}^{m'}$.\par\vspace{1mm}
{\rm (iii)} 	$\htv\,\cyq\,\htv'$ and $m=m'$.
\endPR
	Since (i)$\to$(ii) clearly holds, it is sufficient to show that (iii)$\to$(i) and (ii)$\to$(iii). 
\pvan
(iii) $\to$ (i): 
	By condition that $\htv\,\cyq\,\htv'$ and Claim \ref{clm:flype}, we have that $P_{\htv}^m\approx P_{\htv'}^{m}$.
	Then by condition that $m=m'$, we have that $P_{\htv'}^{m}\approx P_{\htv'}^{m'}$.
\pvan
(ii) $\to$ (iii): 
	By condition (ii) and Claim \ref{clm:sde} (1), each corresponding sublinks of $P_{\htv}^m$ and 
	$P_{\htv'}^{m'}$ are self $\De$-equivalent. The $2$-component non-trivial sublinks of $P_{\htv}^m$ 
	(resp. $P_{\htv'}^{m'}$) are $\tau(k_1\ep_1)$, $\dots$, $\tau(k_\m\ep_\m)$ (resp. $\tau(\ell_1\d_1)$, 
	$\dots$, $\tau(\ell_\m\d_\m)$), since $\m\geq 3$. Hence there is a certain $h$ $(1\leq h\leq\m)$ such that 
	$\tau(\ell_i\ep'_i) = \tau(k_{i+h}\ep_{i+h})$ or $\tau(\ell_i\ep'_i) = \tau(k_{\m-i+h}\ep_{\m-i+h})$ 
	$(i = 1, 2, \dots, \m)$ by Claim \ref{clm:sde} (1) and Corollary \ref{cor:TorusDelta}, since none of $k_i$ and 
	$l_i$ is $0$ nor $-2$ $(i=1,$ $\cdots,$ $\m)$.
	Therefore we have that $\htv \cyq\,\htv'$ and thus, the two integers $m$ and $m'$ have the same parity.
	Now it is sufficient to show that $m=m'$ under the assumption $\htv=\htv'$.
\pva
	Let $P_{\htv}^0$ be the link obtained from $P_{\htv}^{m}$ by smoothing a crossing of parameter $e_i$. 
	Note that $P_{\htv}^0$ is a certain connected sum of $\tau(k_1\ep_1),$ $\dots, \tau(k_\m\ep_\m)$. 
	We have the following by Axiom of Conway polynomial (2), by Corollary \ref{cor:TorusAlex}, and by that
	$\Na_{\l\sharp\l'}=\Na_\l\cdot\Na_{\l'}$, where $\l\sharp\l'$ is a certain connected sum of $\l$ and $\l'$:
	$$a_{\m+1}(P_{\htv}^{m}) -a_{\m+1}(P_{\htv}^{m-2}) 
	= a_{\m}(P_{\htv}^0) 	= \pm {{k_1} \over 2}\times \cdots \times {{k_\m} \over 2}\neq 0.$$
 	Applying this process repeatedly, we have that
	$$\vert a_{\m+1}(P_{\htv}^m) - a_{\m+1}(P_{\htv'}^{m'}) \vert =
  	{{\vert m - m' \vert} \over 2} \times {\vert a_{\m}(P_{\htv}^0) \vert}.$$
	Since any components of $P_{\htv}^m$ and $P_{\htv'}^{m'}$ are trivial, their values of $a_2$ are both $0$.
 	Hence we have that $a_{\m+1}(P_{\htv}^m) = a_{\m+1}(P_{\htv'}^{m'})$ by Proposition \ref{prp:Nak},
 	and thus that $m=m'$.
\endR
%
%
%
%
%
 \bgnR
 	[Proof of Theorem \ref{thm:EPS}.]
	Let $\htv$ and $\htv_e$ be $(k_1\ep_1, k_2\ep_2, \dots, k_u\ep_u)$ and 
	$(k_{i_1}\ep_{i_1},$ $k_{i_2}\ep_{i_2}$, $\dots$, $k_{i_\m}\ep_{i_\m})$, respectively. 
	Moreover, let $\r_o(\htv)=\sum_{k_i: {\rm odd}} k_i$ and $\r_g(\htv)=\sharp \{ i \vert k_i = -2\}$,
	i.e., $\r(\htv)=\r_o(\htv)-\r_g(\htv)$. Define $\htv_s=(k'_{i_1}\ep_{i_1},$ $k'_{i_2}\ep_{i_2}$, 
	$\dots$, $k'_{i_\m}\ep_{i_\m})$ by $k'_{i_j}=k_{i_j}$ if $k_{i_j}\neq -2$ and $k'_{i_j}=-k_{i_j}$ 
	if $k_{i_j}=-2$ $(j=1, 2, \dots, \m)$. Then we have the following:
\par\noindent
	$P_{\htv}$ $\sdq$ $P_{\htv_e}^{\r_o(\htv)}$ 
	$=P(k_{i_1}\ep_{i_1},$ $\dots$, $k_{i_\m}\ep_{i_\m}$, $e_1$, $\dots$, $e_{|\r_o(\htv)|})$ 
	(by Lemma \ref{lem:key}) \pva\hspace{.5mm}
	$\approx$ $P(k'_{i_1}\ep_{i_1},$ $\dots$, $k'_{i_\m}\ep_{i_\m}$, $e_1$, $\dots$, $e_{|\r_o(\htv)|}$, 
	$e'_1$, $\dots$, $e'_{\r_g(\htv)})$, where $e'_i=-1$ (by Claim \ref{clm:normal}) \pva\hspace{.5mm}
	$\approx$ $P_{\htv_s}^{\r(\htv)}=P(k'_{i_1}\ep_{i_1},$ $\dots$, $k'_{i_\m}\ep_{i_\m}$, $e_1$, 
	$\dots$, $e_{|\r(\htv)|})$ (by Reidemeister II moves).
\pva
	Similarly defining $\htv'$ and $\htv'_s$, we have that $P_{\htv'}$ $\sdq$ $P_{\htv'_s}^{\r(\htv')}$.
	Hence $P_{\htv}$ $\sdq$ $P_{\htv'}$ if and only if $P_{\htv_s}^{\r(\htv)}$ $\sdq$ $P_{\htv'_s}^{\r(\htv')}$.
	By Proposition \ref{prp:EPS}, we have that $P_{\htv_s}^{\r(\htv)}$ $\sdq$ $P_{\htv'_s}^{\r(\htv')}$
	if and only if $\htv_s\,\cyq\,\htv'_s$ and  $\r(\htv)=\r(\htv')$. Therefore we know that the statement holds, 
	since $\htv_e=\htv_s$ and $\htv'_e=\htv'_s$ by permitting $2\rr = -2\rs$ and $2\rs = -2\rr$.
\endR
 
%
%
%
%
%
%
%
%
%
%
%
%
%
%
%
%
%
%
%
%

\section{$2$-component pretzel links}\label{sec:2comp}
	In this section, we consider the self $\De$-equivalence for $2$-component pretzel links. By Claim 
	\ref{clm:comp}, a $2$-component pretzel link $\k$ has a diagram $P(k_1\ep_1, k_2\ep_2, \ldots, k_u\ep_u)$ 
	such that all the parameters $k_i$ are odd (odd pretzel link) or exactly two parameters are even (even 
	pretzel link). Then, each of the two knot components is the trivial knot, and thus the value of $a_2$
	of each component is zero. Hence by Proposition \ref{prp:NO}, the self $\De$-equivalence for 
	$2$-component pretzel links are determined by the values of $a_1$ and $a_3$ of the links. 
	Note that we can calculate the Conway polynomials of pretzel links using Lemma \ref{lem:TwRed}.
%
\bgnP \label{prp:CwyPrz}
	Let $\htv=(k_1\ep_1$, $k_2\ep_2$, $\dots$, $k_\m\ep_\m)$. Then we have the following:
	$$\Na_{P_{\htv}}=\Sum\prod_{j=1}^\m \Ga_j(l_j)\Na_{P(l_1\ep_1,l_2\ep_2,\dots, l_\m\ep_\m)},$$
	where the sequence of indices in the sum runs all the $2^\m$ possibilities,  i.e., each $l_i$ takes
	a value $0$ $($resp. $1)$ or $\infty$ if $k_i$ is even $($resp. odd$)$ and $\ep_j=\rs$, $1$ or $0$
	if $\ep_j=\rr$, and $\Ga_j(l_j)=
\bgnseq
1 \mathrm{\, (resp. \,}  -p_jz  \mathrm{)~~~\, if\,} k_j=2p_j,~~~ \ep_j=\rs, l_j=0 \mathrm{\, (resp. \,} \infty \mathrm{)}	 \\[1ex]
1 \mathrm{\, (resp. \,}  -p_jz  \mathrm{)~~~\, if\,} k_j=2p_j+1, \ep_j=\rs, l_j=1\mathrm{\, (resp. \,} \infty \mathrm{)}	\\[1ex]	
\Phi(p_j,z)  \mathrm{\, (resp. \,} \Psi(p_j-1,z)  \mathrm{) ~~~\, if\,} k_j=2p_j, \ep_j=\rr, l_j=1 \mathrm{\, (resp. \,} 0 \mathrm{)} \\[1ex]
\Psi(p_j,z) \mathrm{\, (resp. \,} \Phi(p_j,z)  \mathrm{)~~~\, if\,} k_j=2p_j+1, \ep_j=\rr, l_j=1 \mathrm{\, (resp. \,} 0 \mathrm{)}	 
\endseq$
\endP

\subsection{$2$-component odd pretzel links}
	In this case, $u$ is even, and all the $\ep_i$'s. are $\rs$ or $\rr$. Then, we have the following.
%
\bgnP 
\label	{prp:odd}
		
		Let $L$ be a $2$-compont odd pretzel link which has a diagram $P(k_1\ep, k_2\ep, \ldots, k_u\ep)$.
		Then we have the following, where $k_j=2p_j+1$ and $u=2\n$.   \\
$(1)$		
		If $\ep=\rs$, then $a_1(L) = -\n-\ovs{u}{\uds{j=1}{\sum}} p_j$ and \\
		$a_3(L) = -\lpa \Frc{\n(\n^2-1)}{6}+\Frc{\n(\n-1)}{2}\Sum_{i=1}^{u} p_i
		+(\n-1)\Sum_{1\leq i<j \leq u} p_ip_j+\Sum_{1\leq i<j<k \leq u} p_ip_jp_k\rpa$ 	\pag
		$= \Frc{1}{\,2\,}\,a_1(L)\,  \n(\n-1) + \Frc{1}{6}\,\n(\n-1)(2\n-1)
		-(\n-1)\!\!\!\Sum_{1\leq i<j \leq u} \!\!\! p_ip_j-\!\!\!\Sum_{1\leq i<j<k \leq u} \!\!\! p_ip_jp_k$. 	
		\\
$(2)$		If $\ep=\rr$, then $a_1(L) = \n+\ovs{u}{\uds{j=1}{\sum}} p_j$ and \\
		$a_3(L) ={\Frc{1}{\,2\,}} \,a_1(L) \Sum_{j=1}^u p_j(p_j+1)
			-\Frc{1}{\,6\,} \Sum_{j=1}^u p_j(p_j+1)(2p_j+1)$.
\endPR
	(1) 	We deform all the crossings for $k_i$ into $1$ or $\infty$ by using Lemma \ref{lem:TwRed}. Then, 
		$\Na_{L}$ is represented by a linear combination of $\Na_{\mL}$ and $\Na_{\mL_{t_1,t_2,\ldots,t_\a}}$
		as follows, where $\mL$ is a pretzel link which has a diagram $P(k'_1\rs, k'_2\rs, \ldots, k'_u\rs)$ with 
		$k'_j=1$ for any $j$ $(i=1,\ldots, u)$ and $\mL_{t_1,t_2,\ldots,t_\a}$ is a pretzel link which has a diagram 
		$P(k'_1\rs, k'_2\rs, \ldots, k'_u\rs)$ with $k'_{t_i}=\infty$ and $k'_{j}=1$ if $j\neq t_i$ $(i=1,\ldots, \a)$. 
		Therefore we obtain the conclusion by Proposition \ref{prp:CwyPrz} and Corollary \ref{cor:TorusAlex}.
\pva
		$\Na_{L}=\Na_{\mL}	-\ovs{u}{\uds{i=1}{\sum}}\, p_iz \Na_{\mL_i}
						+\uds{1\leq i<j \leq u}{\sum}p_ip_jz^2 \Na_{\mL_{i,j}}
						-\uds{1\leq i<j<k \leq u}{\sum} p_ip_jp_kz^3 \Na_{\mL_{i,j,k}}+f_1(z)z^4$\pvan
		$=\Na_{\tau(-2\n\rr)}	-\ovs{u}{\uds{i=1}{\sum}}\, p_iz \Na_{\tau(-2\n+1\rr)}
						+\!\!\!	\uds{1\leq i<j \leq u}{\sum} p_ip_jz^2 \Na_{\tau(-2\n+2\rr)}
						-\!\!\!	\uds{1\leq i<j<k \leq u}{\sum}\!\!\!	 
							p_ip_jp_kz^3 \Na_{\tau(-2\n+3\rr)}+f_1(z)z^4$\pvan
		$=-\lpa \n z+\Frc{1}{\,6\,}\,\n(\n^2-1)z^3\rpa-\ovs{u}{\uds{i=1}{\sum}}\, p_iz \lpa 1+\Frc{1}{\,2\,}\n(\n-1)z^2\rpa
					-\uds{1\leq i<j \leq u}{\sum} p_ip_jz^2  (\n-1)z$ \pac
		$\!\!\!			-\uds{1\leq i<j<k \leq u}{\sum} \!\!\!	p_ip_jp_kz^3+f_2(z)z^4$. 
\npg
	(2) 	We deform all the crossings for $k_i$ into $1$ or $0$ by using Lemma \ref{lem:TwRed}. Hence, 
		$\Na_{L}$ is represented by a linear combination of $\Na_{\mL}$ and $\Na_{\mL_{t_1,t_2,\ldots,t_\a}}$
		as follows, where $\mL$ is a pretzel link which has a diagram $P(k'_1\rs, k'_2\rs, \ldots, k'_u\rs)$ with 
		$k'_j=1$ for any $j$ $(i=1,\ldots, u)$ and $\mL_{t_1,t_2,\ldots,t_\a}$ is a pretzel link which has a diagram 
		$P(k'_1\rs, k'_2\rs, \ldots, k'_u\rs)$ with $k'_{t_i}=0$ and $k'_{j}=1$ if $j\neq t_i$ $(i=1,\ldots, \a)$. 
		Therefore we obtain the conclusion by Proposition \ref{prp:CwyPrz} and Corollary \ref{cor:TorusAlex}. 
		Here note that $\mL_{t_1,t_2,\ldots,t_\a}$ is the $\a$ component trivial link, and thus 
		$\Na_{\mL_{t_1}}=1$ and $\Na_{\mL_{t_1,t_2,\ldots,t_\a}}=0$ if $\a\geq 2$. 
\pva
		$\Na_{L}=\ovs{u}{\uds{i=1}{\prod}}\Psi(p_i)\Na_{\mL}
		+ \ovs{u}{\uds{j=1}{\sum}}\, \uds{i\neq j}{\prod}\Psi(p_i)\Phi(p_j)\Na_{\mL_j}$
		$=\ovs{u}{\uds{i=1}{\prod}}\Psi(p_i)\Na_{\tau(-2\n\rs)}
		+ \ovs{u}{\uds{j=1}{\sum}}\, \uds{i\neq j}{\prod}\Psi(p_i)\Phi(p_j)\cdot 1$\pvan
		$=\n z \ovs{u}{\uds{i=1}{\prod}} \lpa 1+\y_1(p_i)z^2\rpa
		 +\ovs{u}{\uds{j=1}{\sum}}\, \uds{i\neq j}{\prod}\lpa 1+\y_1(p_i)z^2\rpa
			\lpa p_j z+\v_1(p_j)z^3\rpa+g_1(z)z^4$\pvan
		$=\n z \lpa 1+\ovs{u}{\uds{i=1}{\sum}} \y_1(p_i)z^2\rpa
		 +\ovs{u}{\uds{j=1}{\sum}}\, \lpa 1+\lpa \ovs{u}{\uds{i=1}{\sum}}\, \y_1(p_i) -\y_1(p_j) \rpa z^2\rpa
			\lpa p_j z+\v_1(p_j)z^3\rpa+g_2(z)z^4$\pvan
	 	$=\n z \lpa 1+\ovs{u}{\uds{i=1}{\sum}} \y_1(p_i)z^2\rpa	+\ovs{u}{\uds{j=1}{\sum}}\, 
		\lpa p_j z+ \lpa p_j \lpa \ovs{u}{\uds{i=1}{\sum}}\, \y_1(p_i) -\y_1(p_j) \rpa
		 +\v_1(p_j)\rpa z^3\rpa	+g_3(z)z^4$\pvan
	 	$=\lpa \n+ \ovs{u}{\uds{j=1}{\sum}}\, p_j\rpa z 
		 + \lpa \ovs{u}{\uds{i=1}{\sum}}\, \y_1(p_i) \lpa \n+ \ovs{u}{\uds{j=1}{\sum}}\, p_j\rpa
		 + \ovs{u}{\uds{j=1}{\sum}}\, \lpa \v_1(p_j) - p_j \y_1(p_j)\rpa	 \rpa z^3+g_3(z)z^4$\pvan
		$=\lpa \n+ \ovs{u}{\uds{j=1}{\sum}}\, p_j\rpa z 
		 + \lpa  \Frc{1}{\,2\,}\ovs{u}{\uds{i=1}{\sum}}\,p_j(p_j+1)
			  \lpa \n+ \ovs{u}{\uds{j=1}{\sum}}\, p_j\rpa
		 - \Frc{1}{\,6\,}\ovs{u}{\uds{j=1}{\sum}}\, p_j(p_j+1)(2p_j+1)	 \rpa z^3+g_3(z)z^4$
\pvc
		The last equation is obtained by the following:
		$$\v_1(p_j) - p_j \y_1(p_j)= \Frc{1}{\,6\,} p_j (p_j^2-1) -  \Frc{1}{\,2\,} p_j^2 (p_j+1)
		=- \Frc{1}{\,6\,} p_j (2p_j^2+3p_j+1).$$
\endR
%
%
%
%
%
\subsection{$2$-component even pretzel links}
		In this case, we may assume that $k_1$ and $k_2$ are even and $(\ep_1,\ep_2)=(\rs,\rs)$, 
		$(\rr,\rr)$, or $(\rs,\rr)$ by Claim \ref{clm:flype}. Then, we have the following.
\bgnP 
\label	{prp:even}
		
		Let $L$ be a $2$-compont even pretzel link which has a diagram 
		$P_{\htv}=P(2p\ep_1, 2q\ep_2, k_3\rr, \ldots, k_u\rr)$. Then we have the following, where 
		$k_j=2p_j+1$, $u=2\n$, and $m=\ovs{u}{\uds{i=3}{\sum}} k_i$.  
$(1)$		If $\ep_1	=\ep_2=\rs$, then $a_1(L) = -(p+q)$ and $a_3(L)= \Frc{1}{\,2\,}mpq$\\
$(2)$		If $\ep_1	=\ep_2=\rr$, then $a_1(L) =   p+q$ and $a_3(L)
				= \Frc{1}{\,6\,}(p+q+1)(p+q)(p+q-1)+\Frc{1}{\,2\,}mpq$\\
$(3)$		If $\ep_1	=\rs$ and $\ep_2=\rr$, then $a_1(L) =  q-p$ and 
				$a_3(L)= \Frc{1}{\,6\,}q(q^2-1)-\Frc{1}{\,2\,}(m+q)pq$
\endPR
		By Lemma \ref{lem:key}, we have that $P_{\htv} \sdm P(2p\ep_1, 2q\ep_2, e_1\rr, \ldots, 
		e_{|m|}\rr)$, where $e_i=\Frc{m}{|m|}$ $(i=3, \ldots, u)$. 
		Let $L(\a_1\ep_1, \b_2\ep_2)=P(\a_1\ep_1, \b_2\ep_2, e_1\rr, \ldots, e_{|m|}\rr)$.
\pvan
	(1) 	We deform all the crossings for $p$ and $q$ into $0$ or $\infty$ by using Lemma \ref{lem:TwRed}. 
		Note that $m$ is even. Since $\l(0\rs,0\rs)$ is the $2$-component trivial link, $\l(0\rs,\infty\rs)$ and 
		$\l(\infty\rs,0\rs)$ are the trivial knots, and $\l(\infty\rs,\infty\rs)=\tau(-m\rs)$, we have the following by 
		Proposition \ref{prp:CwyPrz}, and hence we have the conclusion, since $\Na_{\tau(-m\rs)}=\Frc{1}{2}mz$:
		\pvan
		$\Na_{L}=1\cdot 1 \Na_{\l(0\rs,0\rs)}+(-pz)\cdot 1 \Na_{\l(0\rs,\infty\rs)}
		+1\cdot (-qz)\Na_{\l(\infty\rs,0\rs)}+(-pz)(-qz)\Na_{\l(\infty\rs,\infty\rs)}$\pac
		$=0-(p+q)z+pqz^2\Na_{\tau(-m\rs)}$
\pvan
	(2) 	We deform all the crossings for $p$ and $q$ into $1$ or $0$ by using Lemma \ref{lem:TwRed}. 
		Note that $m$ is even. Since $\l(0\rr,0\rr)$ is the $2$-component trivial link, $\l(1\rr,0\rr)$ and 
		$\l(0\rr,1\rr)$ are the trivial knots, and $\l(1\rr,1\rr)=\tau(-(m+2)\rs)$, we have the following by Proposition 
		\ref{prp:CwyPrz}, and hence we have the conclusion, since $\Na_{\tau(-(m+2)\rs)}=\Frc{1}{2}(m+2)z$.
		\pvan
		$\Na_{L}=\Psi(p-1,z)\Psi(q-1,z)\Na_{\l(0\rr,0\rr)}+\Phi(p,z)\cdot \Psi(q-1,z) \Na_{\l(1\rr,0\rr)}$\pag\hsp{5.4cm}
		$+\,\Psi(p-1,z)\cdot \Phi(q,z)\Na_{\l(0\rr,1\rr)}+\Phi(p,z)\cdot \Phi(q,z) \Na_{\l(1\rr,1\rr)}$\par\hsc
		$=0+(pz+\Frc{1}{\,6\,}p(p^2-1)z^3+\cdots)(	1+\Frc{1}{\,2\,}(q-1)qz^2+\cdots)$\pag\hsd
		$+(1+\Frc{1}{\,2\,}(p-1)pz^2+\cdots)(qz+\Frc{1}{\,6\,}q(q^2-1)z^3+\cdots)$\pag\hsd
		$+(pz+\Frc{1}{\,6\,}p(p^2-1)z^3+\cdots)(qz+\Frc{1}{\,6\,}q(q^2-1)z^3+\cdots)\Frc{1}{\,2\,}(m+2)z$\par\nid
		$=\!(p+q)z+\Frc{1}{\,6\,}\!\lpa 3pq(q-1)\!+\!p(p^2-1)\!+\!3pq(p-1)\!+\!q(q^2-1)+\!6pq\rpa\!z^3\!
			+\Frc{1}{\,2\,}mpqz^3+f(z)z^4$\pac
		$=(p+q)z+\Frc{1}{\,6\,}\lpa 3pq(p+q)+(p^3+q^3)-(p+q)\rpa z^3+\Frc{1}{\,2\,}mpqz^3+f(z)z^4$\pac
		$=(p+q)z+\Frc{1}{\,6\,}(p+q) \lpa (p+q)^2-1\rpa z^3+\Frc{1}{\,2\,}mpqz^3+f(z)z^4$.
\pvan
	(3) 	We deform all the crossings for $p$ and $q$ into $0$, $1$, or $\infty$ by using Lemma \ref{lem:TwRed}. 
		Note that $m$ is odd.
		Since $\l(0\rs,0\rr)$ is the $2$-component trivial link, $\l(0\rs,1\rr)$ and $\l(\infty\rs,0\rr)$ are the 
		trivial knots, and $\l(\infty\rs,1\rr)=T(-(m+1)\rs)$. Thus, we have the following by Proposition 
		\ref{prp:CwyPrz}, and hence we have the conclusion, since $\Na_{T(-(m+1)\rs)}=\Frc{1}{\,2\,}(m+1)z$.
		\pvan
		$\Na_{L}=1\cdot\Psi(q-1,z)\Na_{\l(0\rs,0\rr)}+1\cdot\Phi(q,z)\Na_{\l(0\rs,1\rr)}
		+(-pz)\cdot\Psi(q-1,z)\Na_{\l(\infty\rs,0\rr)}+(-pz)\Phi(q,z)\Na_{\l(\infty\rs,1\rr)}$\pac
		$=0+(qz+\Frc{1}{\,6\,}q(q^2-1)z^3+\cdots)+(-pz)(1+\Frc{1}{\,2\,}(q-1)qz^2+\cdots)$\pag\hsp{5.4cm}
		$+(-pz)(qz+\Frc{1}{\,6\,}q(q^2-1)z^3+\cdots)\Frc{1}{\,2\,}(m+1)z$\pac
		$=(q-p)z+\Frc{1}{\,6\,}q(q^2-1)z^3-\Frc{1}{\,2\,}(m+q)pqz^3$.
\endR

		For example, we consider two $2$-component pretzel links $\k$ and $\l$ with diagrams
		$P(6\rr, -6\rr, 1\rr, 1\rr)$ and $P(4\rs, 4\rr, 1\rr, 1\rr, 1\rr)$, respectively. We have 
		$\Na_\k(z) = -9z^3 - 24z^5 -22z^7 -8z^9 - z^{11}$ and $\Na_\l(z) = -9z^3 -4z^5$.
		Then we have $\k\sdq\l$ by Proposition \ref{prp:NO}.

%
%
%
%
%
%
%
%
%
%
%
%
%
%
%
%
%
%
%
%

\section{Self Delta-equivalence of slice pretzel links}\label{sec:eras}

	In our convention, the product ${\mathbf R}^3 \times (-\infty, 0]$ is oriented so that the natural 
	projections $p_0: {\mathbf R}^3  \times \{0\} \rightarrow {\mathbf R}^3$ is orientation-preserving. 
	Further, when a link $\k$ is given in ${\mathbf R}^3$, the orientations of 
	$\k \times \{0\} (\subset {\mathbf R}^3  \times \{0\})$ is specified so that 
	$p_0 \vert_{\k \times \{0\}}: \k \times \{0\} \rightarrow \k$ is orientation-preserving. 
	Thus, the link  $({\mathbf R}^3  \times \{0\}, \k \times \{0\})$ is identified with $({\mathbf R}^3 , \k)$.
	The link $\k$ is {\it slice\/} if there is a locally flat, oriented, and proper $\m$ disks $C$ such that 
	$C \cap {\mathbf R}^3  \times \{0\} = \k \times \{0\}$.
\pva           
	For an even number sequence of non-zero integers $k_i$ $(i = 1, 2, \dots, 2v)$,
	we define an {\it erasable sequence} $\vv= (k_1,k_2, \dots, k_{2v})$ inductively as follows.
(1) 	If $m = 1$, then $k_1 + k_2 = 0$.
(2) 	For $v \geq 2$,
	there are two integers $k_i$ and $k_j$  $(1 \leq i < j \leq 2v)$ in $\vv$ such that $k_i + k_j = 0$ 
	and the sequence with $2v-2$ elements obtained from $\vv$ by deleting $k_i$ and $k_j$ is erasable. 
	The following is known.

\bgnP{\rm(\cite[Theorem~1.3]{KST20}, \cite[Corollary~1.5]{NST})}\label{prop27}
		Let $\k$ be a slice pretzel link. Then $\k$ is equivalent to one of the following.

{\rm (i)} 	A slice knot.

{\rm (ii)} 	A $2$-component pretzel link with a diagram $P(k_1\ep_1, k_2\ep_2, \ldots, k_u\ep_u)$ 
		such that $u$ is even, $(k_1, k_2, \dots, k_{u})$ is erasable, and $|k_i|\geq 2$ for each 
		$i$ $(i,=1,2,\ldots, u)$.

\endP

\bgnR
	[Proof of Theorem \ref{thm:Slice}.]  \label{sec:thm14}
	By Proposition \ref{prop27}, we have two cases (i) and (ii). In Case (i), we know that any (slice) knot 
	is self $\De$-equivalent to the trivial knot. In Case (ii), the slice pretzel link $\k$ has a $2$ components 
	and a diagram $P(k_1\ep_1,$ $k_2\ep_2,$ $\ldots,$ $k_{2v}\ep_{2v})$ such that $(k_1,$ $k_2,$ $\dots,$ 
	$k_{2v})$ is erasable, and $|k_i|\geq 2$ for each $i$ $(i,=1,2,\ldots, u)$ by Proposition \ref{prop27}.
	Hence, the linking number of $\k$, i.e., $a_1(\k)$ is $0$, and thus $a_3(\k)$ is a concordance invariant
	(\cite[Theorem~3.2]{Coc85}). Since $\k$ is slice, $a_3(\k)=0$. Therefore, $\k$ is self $\De$-equivalent 
	to the $2$-component trivial link by Proposition \ref{prp:Nak}.
\endR

	In the above proof, it is in fact the case that $\Na_\k(z)=0$, since the first non-vanishing coefficient of the 
	Conway polynomial is a concordance invariant (\cite[Theorem~3.2]{Coc85}). Furthermore, the condition 
	that $(k_1,$ $k_2,$ $\dots,$ $k_{2v})$ is erasable is sufficient to imply $\Na_\k(z)=0$.
	
\bgnP	\label{prp:spconway}
	Let $\k$ be a $2$-component link with a diagram $P_{\vv}$ such that $\vv=(k_1,$ $k_2,$ $\dots,$ $k_{2v})$ 
	is erasable. Then we have that $\Na_\k(z)=0$.
\endP
		
	For two link diagrams $K=P(\dots, k_i\ep_i, k_{i+1}\ep_{i+1}, \dots)$ and $L=P(\dots, k_{i+1}\ep_{i+1}, 
	k_{i}\ep_{i}, \dots)$ as in Figure \ref{fig:mutation}, a local move between $K$ and $L$ is one of mutations.
	Here note that each $\ep_i$ does not change by mutations, since we reverse the orientations of all the 
	components inside the dotted circle illustrated in Figure \ref{fig:mutation} if the orientations of strands 
	inside and outside of the dotted circle do not agree. It is known that any mutation never changes the 
	Conway polynomials \cite{L}.

\bgnF
		\iclg{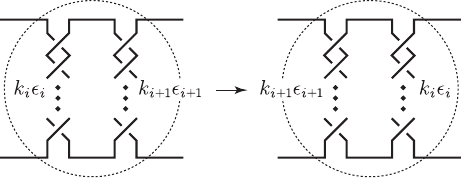} 	\caption{mutation} 	\label{fig:mutation}
\endF
%
\bgnP
\label	{prp:mutation}
		Any mutation on a link $\ell$ does not change the Conway polynomial of $\ell$.
\endP

	Let $\htv=(k_1\ep_1, k_2\ep_2, \ldots, k_{2v}\ep_{2v})$ be an enhanced pretzel sequence such that 
	$\vv= (k_1, k_2,$ $\dots,$ $k_{2v})$ is erasable and $|k_i|\geq 2$ $(i=1,$ $2,$ 	$\ldots,$ $2v$).
	Since $\vv$ is erasable, there exists a permutation $\s$ which is a composition of $v$ 
	mutually disjoint transpositions such that $k_i+k_{\s(i)}=0$ for any $i$ $(i=1, 2, \cdots, 2v)$. 
	We call $\s$ a {\it pairing permutation for} $\vv$.

\bgnP
\label	{prp:TrvCwy}
		Let $\htv=(k_1\ep_1, k_2\ep_2, \ldots, k_{2v}\ep_{2v})$ be an enhanced pretzel sequence such that 
		$\vv= (k_1, k_2,$ $\dots,$ $k_{2v})$ is erasable and $|k_i|\geq 2$ $(i=1,$ $2,$ 	$\ldots,$ $2v)$.
		If a pairing permutation $\s$ for $\vv$ satisfies that $\ep_{\s(i)} = \ep_{i}$ for each $i $ $(i=1,$ $2,$ 
		$\ldots,$ $2v)$, then we have that $\Na_{\k}(z) = 0$, where $\k=P_{\htv}$.
\endP
	The condition that $\ep_{\s(i)} = \ep_{i}$ is necessary. For instance, if $\htv=(2\rs,$ $-2\rr,$
	$2\rs,$ $-2\rr)$, then $\Na_{\k}(z)=-4z^3-z^5\neq 0$, where $\k=P_{\htv}$.
\bgnR
	By the given condition, we can reorder $\htv$ into
	$\htv'= (k_{i_1}\ep_{i_1}$, $-k_{i_1}\ep_{i_1}$, $\ldots$, $k_{i_{v}}\ep_{i_{v}}$, $-k_{i_{v}}\ep_{i_{v}})$.
	By Proposition \ref{prp:mutation}, it is sufficient to show that $\Na_{\l}(z)=0$, where $\l=P_{\htv'}$.
	We show the statement by induction on $v$. 
	If $v=1$, then we have that $P(k_{i_1}\ep_{i_1}$, $-k_{i_1}\ep_{i_1})=\tau(2,k_{i_1}-k_{i_1})$ is 
	the $2$-component trivial link $O^2$, and thus that $\Na_{\l}(z)=0$.
\pva
	Suppose that $v\geq 2$ and that the statement holds for $v-1$. By Proposition \ref{prp:ErsCwy},
	$\Na_{\l}(z)$ is the product of $\Na_{\l'}(z)$ and $-p_{i_v}^2z^2$, $(1+p_{i_v}^2z^2)$, 
	$-\Phi^2(p_{i_v},z)$, or $\Psi^2(p_{i_v},z)$ depending on the parities of $k_{i_v}$ ($=2p_{i_v}+1$) 
	and $\ep_{i_v}$, where $\l'=P(k_{i_1}\ep_{i_1}$, $-k_{i_1}\ep_{i_1}$, $\ldots$, 
	$k_{i_{v-1}}\ep_{i_{v-1}}$, $-k_{i_{v-1}}\ep_{i_{v-1}})$. Since $\Na_{\l'}(z)=0$ by the assumption, 
	we have that $\Na_{\l}(z)=0$.
\endR	

\bgnR
	[Proof of Proposition \ref{prp:spconway}.]  
	By Proposition \ref{prop27}, $\k$ has a diagram $P(k_1\ep_1,$ $k_2\ep_2,$ $\ldots,$ $k_{2v}\ep_{2v})$ 
	such that $(k_1,$ $k_2,$ $\dots,$ $k_{2v})$ is erasable, and $|k_i|\geq 2$ for each $i$ $(i,=1,2,\ldots, u)$.
	Then $\htv$ has no even parameters or exactly two even parameters by 
	Claim \ref{clm:comp}. In the former case, $\ep_i$'s are all $\rs$ or all $\rr$ by Proposition \ref{prp:ori} (2).
	In the latter case, let $k_a$ and $k_b$ be the two even parameters.	Then $\ep_a$ and $\ep_b$ are both 
	$\rs$ or both $\rr$ and the other $2(v-1)$ $\ep_i$'s are all $\rr$ by Proposition \ref{prp:ori} (3). Thus in 
	both cases, any pairing permutation $\s$ for $\htv$ satisfies that $\ep_{\s(i)} = \ep_{i}$ for each $i $ 
	$(i=1,$ $2,$ $\ldots,$ $2v)$. Therefore we have that $\Na_{\k}(z) = 0$ by Proposition \ref{prp:TrvCwy}.
\endR

%
%
%
%
%
%
%
%
%
%
%
%
%
%
%
%
%
%
%
%

\pvg
\noindent{Yasutaka NAKANISHI}{\ }(nakanisi@math.kobe-u.ac.jp)\\
{Department of Mathematics, Kobe University, 
Rokkodai-cho 1-1, Nada-ku, Kobe 657-8501, Japan}

\pvc

\noindent{Tetsuo SHIBUYA}\\
{Department of Mathematics, Osaka Institute of Technology, 
5-16-1 Omiya, Asahi-ku, Osaka 535-8585, Japan}

\pvc

\noindent{Tatsuya TSUKAMOTO}{\ }(tatsuya.tsukamoto@oit.ac.jp)\\
{Department of Mathematics, Osaka Institute of Technology, 
5-16-1 Omiya, Asahi-ku, Osaka 535-8585, Japan}

\end{document}